\documentclass[12pt, twoside, geqno]{article}
\usepackage{amsmath,amsthm}
\usepackage{amssymb,latexsym}
\usepackage{enumerate}
\usepackage{mathrsfs}
\usepackage{amsmath,amssymb}
\newcommand{\be}{\begin{equation}}
\newcommand{\ee}{\end{equation}}
\newcommand{\bal}{\begin{aligned}}
\newcommand{\eal}{\end{aligned}}
\newcommand{\bee}{\begin{equation*}}
\newcommand{\eee}{\end{equation*}}
\usepackage{color}

\def\cwedge{\bigcirc\kern-1.07em\wedge\ }

\newtheorem{thm}{Theorem}[section]

\newtheorem{cor}[thm]{Corollary}

\newtheorem*{ques}{Question}
\newtheorem*{pf}{Proof}
\newtheorem*{rem}{Remark}
\newtheorem*{ex}{Example}

\newtheorem*{mainthm}{Main Theorem}
\numberwithin{equation}{section}

\begin{document}
\begin{center}
{\large \bf A Curvature identity on a 4-dimensional Riemannian manifold}\\
\end{center}
\footnotetext{\small{\it E-mail addresses}: {\bf
prettyfish@skku.edu} (Y. Euh), {\bf parkj@skku.edu} (J. H. Park),
{\bf sekigawa@math.sc.niigata-u.ac.jp} (K. Sekigawa).}

\begin{center}
Yunhee Euh${}^{\dag}$, JeongHyeong Park ${}^{\ddag}$ and Kouei
Sekigawa$^{\dag}$

\end{center}
{\small
\begin{center}
$~~~{}^{\dag}$Department of Mathematics,
    Niigata University,
    Niigata 950-2181, JAPAN\\
$~~~{}^{\ddag}$ Department of Mathematics,
    Sungkyunkwan University,
    Suwon 440-746, KOREA
\end{center}
}

\begin{abstract}
We give a curvature identity derived from the generalized
Gauss-Bonnet formula for 4-dimensional compact oriented Riemannian
manifolds. We prove that the curvature identity holds on any
4-dimensional Riemannian manifold which is not necessarily compact.
We also provide some applications of the identity.
\end{abstract}
\noindent {\it Mathematics Subsect Classification (2010)} : 53B20, 53C20 
\section{Introduction}

We recall that the Gauss-Bonnet formula for a compact oriented
surface $M=(M,g)$
    \begin{equation}\label{eq:G_B_2}
    2\pi\chi(M)=\int_M K dv_g,
    \end{equation}
where $\chi(M)$ is the Euler number of $M$, $K$ is the Gaussian
curvature of $M$ and $dv_g$ is the volume element of $M$. The
Gaussian curvature $K$ of $M$ is also expressed by using the scalar
curvature $\tau$ of $M$ as $K=\frac{\tau}{2}$.

Now we consider any one-parameter smooth deformation $g(t)$ of $g$.
Since the Euler number $\chi(M)$ is a topological invariant of $M$,
from \eqref{eq:G_B_2}, we have
    \be\label{eq:dif_G_B_2}
    0=\frac{d}{dt}\Big|_{t=0}\int_M \tau(t) dv_{g(t)}=\int_M\Big(-\rho^{ij}+\frac{\tau}{2} g^{ij}\Big)h_{ij} dv_g
    \ee
for symmetric $(0,2)$-tensor field
$h_{ij}=\frac{d}{dt}\Big|_{t=0}g(t)_{ij}$ (for more details, refer
to the next section). Thus we have
    \be\label{eq:eins_2}
    \rho_{ij}=\frac{\tau}{2} g_{ij}.
   \ee
It is well-known that the equality \eqref{eq:eins_2} holds for any
2-dimensional Riemannian manifold without the compactness
assumption.

Motivated by the above observation, the following question will
naturally arise.
\begin{ques}
As mentioned above, does the above phenomenon also occur for any
$2n(n\geq2)$-dimensional Riemannian manifold?
\end{ques}
 Concerning the Question, Berger \cite{Beg}
discussed the Gauss-Bonnet formula for the $4$-dimensional compact
Riemannian manifold $(M,g)$ from the variational theoretic
viewpoint. From one of his results (\cite{Beg}, pp. 292), and taking
account of the well-known fact that the Euler number $\chi(M)$ is a
topological invariant of $M$, we can see that the following
curvature identity holds on any $4$-dimensional compact  Riemannian
manifold $(M,g)$:
   \be\label{main equation}
    \check{R}-2\check{\rho}-L\rho+\tau\rho-\frac{1}{4}(|R|^2-4|\rho|^2+\tau^2)g=0.
    \ee
    Here,
     \bee
     \begin{gathered}
     \check{R}:\check{R}_{ij}=\sum_{a,b,c}R_{abci}R^{abc}_{~~~j},\qquad
     \check{\rho}:\check{\rho}_{ij}=\sum_a\rho_{ai}\rho^{a}_{~j},\\
     L:(L\rho)_{ij}=2\sum_{a,b}R_{iabj}\rho^{ab},
     \end{gathered}
    \eee
    where $R$ is the curvature tensor of $M$ and  $\rho$ is the Ricci
    tensor of $M$.\vskip0.3cm

It is a remarkable fact that the identity is a quadratic equation of
the curvature tensor which does not involve the covariant
derivatives of the curvature tensor. Recently, Labbi \cite{La}
extended the above curvature identity to the higher dimensional
cases by using an elegant method. He considered the only compact
case in his paper, however his equality (10) on p.178 allows us to
deduce that his equation also applies to the non-compact case. One
just have to use a purely algebraic computations in the ring of
double forms. As a final result, we shall prove  the following the
Main Theorem.
\begin{mainthm}
Equation \eqref{main equation} holds on any 4-dimensional Riemannian
manifold.
\end{mainthm}
\begin{rem}
        We can check that equation \eqref{main equation} is valid for any
        4-dimensional pseudo-Riemannian manifold with the aid of
        ``Mathematica".
\end{rem}
 In the present paper, first we shall review
Berger's arguments in \cite{Beg} and also give a direct proof of the
Main Theorem in view of its applications to the related topics.\\

The authors would like to express their thanks to Professor Kowalski
for his helpful comments.

\section{Preliminaries}
In this section, we prepare some fundamental formulas derived from
one-parameter deformations of Riemannain metrics and further
introduce a curvature identity on a 4-dimensional  compact oriented
Riemannian manifold derived from the generalized Gauss-Bonnet
formula by making use of those fundamental formulas.

Let $M=(M,g)$ be an $n$-dimensional Riemannian manifold and
$\mathfrak{X}(M)$ the Lie algebra of all smooth vector fields on
$M$. We denote the Levi-Civita connection, the curvature tensor, the
Ricci tensor and the scalar curvature of $M$ by $\nabla$, $R$,
$\rho$ and $\tau$, respectively. The curvature tensor is defined by
$R(X,Y)Z=[\nabla_X,\nabla_Y]Z-\nabla_{[X,Y]}Z$ for $X$, $Y$,
$Z\in\mathfrak{X}(M)$.

Now, we denote by $\mathfrak{M}(M)$ the space of all Riemannian
metrics on $M$. Let $g(t)\in\mathfrak{M}(M)$ be a smooth curve
through $g(0)=g$. We shall also call it a one-parameter deformation
of $g$. Let $(U;x^1, \cdots, x^n)$ be a local coordinate system on a
coordinate neighborhood $U$ of $M$. With respect to the natural
frame $\Big\{\partial_i=\frac{\partial}{\partial
x^i}\Big\}_{i=1,2,\cdots n}$,  we set
$g(t)(\partial_i,\partial_j)=g(t)_{ij}$,
$R(t)(\partial_i,\partial_j)\partial_k=R(t)_{ijk}^{~~~l}\partial_l$,
$\rho(t)(\partial_i,\partial_j)=\rho(t)_{ij}$,
$\tau(t)=g(t)^{ij}\rho(t)_{ij}$ and $(g(t)^{ij})=(g(t)_{ij})^{-1}$.
In particular, we have $g(0)_{ij}=g_{ij}$,
${R(0)_{ijk}}^l=R_{ijk}^{~~~l}$, $\rho(0)_{ij}=\rho_{ij}$ and
$\tau(0)=\tau$.

 We set
    \be\label{eq:dif_g}
    \frac{d}{dt}\Big|_{t=0}g(t)_{ij}=h_{ij}.
    \ee
Then we see that the $(0,2)$-tensor field $h=(h_{ij})$ is symmetric
on $M$ and we also have
    \be\label{eq:dif_inv_g}
    \frac{d}{dt}\Big|_{t=0}g(t)^{ij}=-h^{ij},
    \ee
where we adopt the standard notational convention of tensor
analysis; for example $h^{ij}=g^{ia}g^{jb}h_{ab}$ and so on. We
denote by $dv_{g(t)}$ the volume element of $(M,g(t))$. Then, we
have
    \be\label{eq:dif_vol}
    \frac{d}{dt}\Big|_{t=0}dv_{g(t)}=\frac12g^{ij}h_{ij}dv_g.
    \ee
From \eqref{eq:dif_g} and \eqref{eq:dif_inv_g}, we see that the
coefficient $\Gamma(t)_{ij}^{~~k}$ of $\nabla^{(t)}$, where
$\nabla^{(t)}$ denotes the Levi-Civita connection with respect to
the metric $g(t)$ satisfies
    \be\label{eq:dif_Cristoffel}
    \frac{d}{dt}\Big|_{t=0}\Gamma(t)_{ij}^{~~k}=\frac12g^{ka}(\nabla_ih_{aj}+\nabla_jh_{ia}-\nabla_ah_{ij}).
    \ee
Thus, from \eqref{eq:dif_Cristoffel}, the derivatives of
$R(t)_{ijk}^{~~~l}$, $\rho(t)_{ij}$ and $\tau(t)$ at $t=0$ are given
respectively by
    \be\label{eq:dif_curv}
    \bal
    \frac{d}{dt}\Big|_{t=0}R(t)_{ijk}^{~~~l}=&\frac12\Big(-R_{ijk}^{~~~a}h_{a}^{~l}+R_{ija}^{~~~l}h_{k}^{~a}\\
    &+\nabla_i\nabla_kh_{j}^{~l}-\nabla_j\nabla_kh_{i}^{~l}-\nabla_i\nabla^lh_{jk}+\nabla_j\nabla^lh_{ik}\Big),
    \eal
    \ee
    \be\label{eq:dif_Ric}
    \bal
    \frac{d}{dt}\Big|_{t=0}\rho(t)_{ij}=&\frac12\Big(-R_{aij}^{~~~b}h_{b}^{~a}+\rho_{ia}h_{j}^{~a}\\
    &+\nabla_a\nabla_{j}h_{i}^{~a}-\nabla_i\nabla_{j}h_{a}^{~a}-\nabla^a\nabla_{a}h_{ij}+\nabla_i\nabla_ah_{j}^{~a}\Big),
    \eal
    \ee
    \be\label{eq:dif_scalar}
    \frac{d}{dt}\Big|_{t=0}\tau(t)=-\rho_{ij}h^{ij}+\nabla^i\nabla^jh_{ij}-\nabla^i\nabla_ih_{j}^{~j}.
    \ee
Equation \eqref{eq:dif_G_B_2} is derived from \eqref{eq:dif_vol} and
\eqref{eq:dif_scalar}. We refer to \cite{LPS} in detail.

In the remainder of this section, we review Berger's result
\cite{Beg}. Let $M=(M,g)$ be a 4-dimensional compact oriented
Riemannian manifold. Then it is well-known that the Euler number
$\chi(M)$ of $M$ is given by the following integral formula
    \be\label{eq:G_B_4}
    \chi(M)=\frac{1}{32\pi^2}\int_M \{|R|^2-4|\rho|^2+\tau^2\}dv_g,
    \ee
where $|R|^2$ and $|\rho|^2$ are the square norms of the curvature
tensor and the Ricci tensor, respectively (namely,
$|R|^2=-g^{ai}g^{bj}R_{ijk}^{~~~l}R_{abl}^{~~~k}$ and
$|\rho|^2=g^{ai}g^{bj}\rho_{ab}\rho_{ij}$). It is well-known as the
generalized Gauss-Bonnet formula. Since $\chi(M)$ is a topological
invariant of $M$, it follows that
    \be\label{eq:dif_G_B_4}
     0=\frac{d}{dt}\Big|_{t=0}\int_M \{ |R(t)|^2-4|\rho(t)|^2+\tau(t)^2\}dv_{g(t)}.
    \ee
Here, from \eqref{eq:dif_inv_g}, \eqref{eq:dif_vol},
\eqref{eq:dif_curv} and taking account of Green's theorem and Ricci
and Bianchi identities, we get
    \be\label{eq:dif_curv_norm}
    \bal
    &\frac{d}{dt}\Big|_{t=0}\int_M|R(t)|^2dv_{g(t)}\\
    =&\int_M\Big\{h^{ai}g^{bj}R_{ijk}^{~~~l}R_{abl}^{~~~k}+g^{ai}h^{bj}R_{ijk}^{~~~l}R_{abl}^{~~~k}-\frac{1}{2}g^{ai}g^{bj}\Big(-R_{ijk}^{~~~c}h_c^{~l}\\
    &\qquad+R_{ijc}^{~~~l}h_k^{~c}+\nabla_i\nabla_kh_j^{~l}-\nabla_j\nabla_kh_i^{~l}-\nabla_i\nabla^lh_{jk}+\nabla_j\nabla^lh_{ik}\Big)R_{abl}^{~~~k}\\
    &\qquad-\frac{1}{2}g^{ai}g^{bj}R_{ijk}^{~~~l}\Big(-R_{abl}^{~~~c}h_c^{~k}+R_{abc}^{~~~k}h_l^{~c}+\nabla_a\nabla_lh_b^{~k}-\nabla_b\nabla_lh_a^{~k}\\
    &\qquad-\nabla_a\nabla^kh_{bl}+\nabla_b\nabla^kh_{al}\Big)+\frac{1}{2}|R|^2g^{ij}h_{ij}\Big\}dv_g.\\
    =&\int_M\Big\{-\Big(2R^{abci}R_{abc}^{~~~j}+4\nabla^a\nabla_a\rho^{ij}-2\nabla^j\nabla^i\tau-4\rho^{i}_{~a}\rho^{ja}+4\rho^{ab}R_{~ab}^{i~~~j}\Big)\\
    &\qquad+\frac{1}{2}|R|^2g^{ij}\Big\}h_{ij}dv_g.
    \eal
    \ee
Similarly, from \eqref{eq:dif_inv_g}, \eqref{eq:dif_vol} and
\eqref{eq:dif_Ric}, we get (by taking account of Green's theorem and
the
 Ricci identity)
    \be\label{eq:dif_rho}
    \bal
    &\frac{d}{dt}\Big|_{t=0}\int_M|\rho(t)|^2dv_{g(t)}\\
    =&\int_M\Big\{-h^{ai}g^{bj}\rho_{ab}\rho_{ij}-g^{ai}h^{bj}\rho_{ab}\rho_{ij}+\frac{1}{2}g^{ai}g^{bj}\Big(-R_{uab}^{~~~v}h_{v}^{~u}+\rho_{au}h_{b}^{~u}\\
    &\qquad+\nabla_u\nabla_bh_{a}^{~u}-\nabla_a\nabla_bh_{u}^{~u}-\nabla^u\nabla_uh_{ab}+\nabla_a\nabla_uh_b^{~u}\Big)\rho_{ij}\\
    &\qquad+\frac{1}{2}g^{ai}g^{bj}\rho_{ab}\Big(-R_{uij}^{~~~v}h_{v}^{~u}+\rho_{iu}h_{j}^{~u}+\nabla_u\nabla_{j}h_{i}^{~u}-\nabla_i\nabla_{j}h_{u}^{~u}\\
    &\qquad-\nabla^u\nabla_{u}h_{ij}+\nabla_i\nabla_uh_{j}^{~u}\Big)+\frac12|\rho|^2g^{ij}h_{ij}\Big\} dv_g\\
    =&\int_M\Big\{-2\rho^{ab}R^{i~~~j}_{~ab}-\frac12(\triangle\tau) g^{ij}
    -\nabla^a\nabla_a\rho^{ij}+\nabla^j\nabla^i\tau+\frac12|\rho|^2g^{ij}\Big\}h_{ij}dv_g,
    \eal
    \ee
where $\triangle$ is the Laplace-Beltrami operator acting on
 differentiable functions on $M$.
Further, from \eqref{eq:dif_inv_g}, \eqref{eq:dif_vol} and
\eqref{eq:dif_scalar}, we get
    \be\label{eq:dif_tau}
    \bal
    &\frac{d}{dt}\Big|_{t=0}\int_M\tau(t)^2 dv_{g(t)}\\
    =&\int_M
    \Big\{2\tau\Big(-\rho_{ij}h^{ij}+\nabla^j\nabla^ih_{ij}-\nabla^i\nabla_ih_{j}^{~j}\Big)+\frac{1}{2}\tau^2g^{ij}h_{ij}
    \Big\}dv_g\\
    =&\int_M\{-2\tau\rho^{ij}+2\nabla^j\nabla^i\tau-2(\triangle\tau)g^{ij}+\frac12\tau^2g^{ij}\}h_{ij}dv_g.
    \eal
    \ee
From \eqref{eq:dif_G_B_4}$\sim$\eqref{eq:dif_tau}, we see that the
following integral formula
    \be\label{eq:dif_G_B_4_id}
    \bal
    &\int_M\Big\{R^{abci}R_{abc}^{~~~j}-2\rho^{ia}\rho_{a}^{~j}-2\rho^{ab}R_{~ab}^{i~~~j}\\
    &\quad+\tau\rho^{ij}-\frac{1}{4}|R|^2g^{ij}+|\rho|^2g^{ij}-\frac{\tau^2}{4}g^{ij}\Big\}h_{ij}dv_g=0
    \eal
    \ee
holds for any symmetric $(0,2)$-tensor field $h=(h_{ij})$.
Therefore, we see that finally the curvature identity
    \be\label{eq:ind}
    \bal
    R^{abci}R_{abc}^{~~~j}-2\rho^{ia}\rho_{a}^{~j}-2\rho^{ab}R_{~ab}^{i~~~j}+\tau\rho^{ij}-\frac{1}{4}\big(|R|^2-4|\rho|^2+\tau^2\big)g^{ij}=0
    \eal
    \ee
holds on $M$. Contracting \eqref{eq:ind} with $g_{iu}g_{jv}$ we can
confirm that Main Theorem is valid for the compact case. Now, let
$\{e_i\}$ be an orthonormal basis of $T_pM$ at any point $p\in M$.
Then, we may rewrite \eqref{eq:ind} as follows:
    \be\label{eq:id}
    \bal
    &\sum_{a,b,c}R_{abci}R_{abcj}-2\sum_{a}\rho_{ia}\rho_{ja}-2\sum_{a,b}\rho_{ab}R_{iabj}\\
    &\qquad+\tau\rho_{ij}-\frac{1}{4}\big(|R|^2-4|\rho|^2+\tau^2\big)\delta_{ij}=0,
    \eal
    \ee
where $R_{ijkl}=g(R(e_i,e_j)e_k, e_l)$, $\rho_{ij}=\rho(e_i,e_j)$.
In this paper, we shall adopt the notational conventions with
respect to a natural basis and an orthonormal basis alternatively
for Main Theorem components of tensors.

From \eqref{eq:dif_curv_norm}, \eqref{eq:dif_rho} and
\eqref{eq:dif_tau}, we see that each term of \eqref{eq:dif_G_B_4}
contains the covariant derivatives of the Ricci tensor and scalar
curvature but equation \eqref{eq:dif_G_B_4} no longer involves the
covariant derivatives as equation \eqref{eq:dif_G_B_4_id}. Based on
this observation Berger inquired whether the similar phenomenon
would hold true for the higher dimension. In \cite{La}, Labbi
recently gave the positive answer to this question.

\section{Proof of Main Theorem}

In this section, we shall give a proof of Main Theorem, namely, the
curvature equality \eqref{eq:id} holds for any (not necessarily
compact) 4-dimensional Riemannian manifold.

Let $M=(M,g)$ be a 4-dimensional Riemannian manifold and $\{e_i\}$ a
Chern basis of $T_pM$ at any point $p\in M$, namely, an orthonormal
basis of $T_pM$ satisfying
    \be\label{eq:Chern_basis}
    R_{1213}=R_{1214}=R_{1223}=R_{1224}=R_{1314}=R_{1323}=0
    \ee
\cite{K-P}. We note that a Singer-Thorpe basis for a 4-dimensional
Einstein manifold is a special kind of a Chern basis \cite{S-T}.
From \eqref{eq:Chern_basis}, we get
    \be\label{eq:curv_1}
    \bal
    \sum_{a,b,c}R_{abc1}^{~~~~2}
    =&2\{R_{1212}^{~~~~2}+R_{1313}^{~~~~2}+R_{1414}^{~~~~2}+R_{1234}^{~~~~2}+R_{1324}^{~~~~2}+R_{1423}^{~~~~2}\\
    &+\rho_{12}^{~~2}+\rho_{13}^{~~2}+\rho_{14}^{~~2}\}.
    \eal
    \ee
Similarly, we get
    \be\label{eq:curv_2}
    \bal
    \sum_{a,b,c}R_{abc2}^{~~~~2}
    =&2\{R_{1212}^{~~~~2}+R_{2323}^{~~~~2}+R_{2424}^{~~~~2}+R_{1234}^{~~~~2}+R_{1324}^{~~~~2}+R_{1423}^{~~~~2}\\
    &+\rho_{12}^{~~2}+\rho_{23}^{~~2}+\rho_{24}^{~~2}+2\rho_{34}^{~~2}\},\\
    \sum_{a,b,c}R_{abc3}^{~~~~2}
    =&2\{R_{1313}^{~~~~2}+R_{2323}^{~~~~2}+R_{3434}^{~~~~2}+R_{1234}^{~~~~2}+R_{1324}^{~~~~2}+R_{1423}^{~~~~2}\\
    &+\rho_{13}^{~~2}+2\rho_{14}^{~~2}+\rho_{23}^{~~2}+2\rho_{24}^{~~2}+\rho_{34}^{~~2}\},\\
    \sum_{a,b,c}R_{abc4}^{~~~~2}
    =&2\{R_{1414}^{~~~~2}+R_{2424}^{~~~~2}+R_{3434}^{~~~~2}+R_{1234}^{~~~~2}+R_{1324}^{~~~~2}+R_{1423}^{~~~~2}\\
    &+2\rho_{12}^{~~2}+2\rho_{13}^{~~2}+\rho_{14}^{~~2}+2\rho_{23}^{~~2}+\rho_{24}^{~~2}+\rho_{34}^{~~2}\}.
    \eal
    \ee
From \eqref{eq:curv_1} and \eqref{eq:curv_2}, we have also
    \be\label{eq:norm_curv}
    \bal
    |R|^2=&4\{R_{1212}^{~~~~2}+R_{1313}^{~~~~2}+R_{1414}^{~~~~2}+R_{2323}^{~~~~2}+R_{2424}^{~~~~2}+R_{3434}^{~~~~2}+2R_{1234}^{~~~~2}\\
    &+2R_{1324}^{~~~~2}+2R_{1423}^{~~~~2}+2\big(\rho_{12}^{~~2}+\rho_{13}^{~~2}+\rho_{14}^{~~2}+\rho_{23}^{~~2}+\rho_{24}^{~~2}+\rho_{34}^{~~2}\big)\}.
    \eal
    \ee
Further, we get the following
    \be\label{eq:Ric}
    \bal
    &\sum_a\rho_{1a}^{~~2}
    =R_{1212}^{~~~~2}+R_{1313}^{~~~~2}+R_{1414}^{~~~~2}+2R_{1212}R_{1313}+2R_{1313}R_{1414}\\
    &\qquad\qquad+2R_{1212}R_{1414}+\rho_{12}^{~~2}+\rho_{13}^{~~2}+\rho_{14}^{~~2},\\
    &\sum_{a,b}\rho_{ab}R_{1ab1}
    =R_{1212}^{~~~~2}+R_{1313}^{~~~~2}+R_{1414}^{~~~~2}+R_{1212}R_{2323}+R_{1212}R_{2424}\\
    &\qquad\qquad\qquad+R_{1313}R_{2323}+R_{1313}R_{3434}+R_{1414}R_{2424}+R_{1414}R_{3434},\\
    &\tau\rho_{11}=2\{R_{1212}^{~~~~2}+R_{1313}^{~~~~2}+R_{1414}^{~~~~2}+2R_{1212}R_{1313}+2R_{1313}R_{1414}\\
    &\qquad+2R_{1212}R_{1414}+R_{1212}R_{2323}+R_{1212}R_{2424}+R_{1212}R_{3434}+R_{1313}R_{2323}\\
    &\qquad+R_{1313}R_{2424}+R_{1313}R_{3434}+R_{1414}R_{2323}+R_{1414}R_{2424}+R_{1414}R_{3434}\},\\
    &|\rho|^2=\rho_{11}^{~~2}+\rho_{22}^{~~2}+\rho_{33}^{~~2}+\rho_{44}^{~~2}+2(\rho_{12}^{~~2}+\rho_{13}^{~~2}+\rho_{14}^{~~2}+\rho_{23}^{~~2}+\rho_{24}^{~~2}+\rho_{34}^{~~2})\\
    &\quad~~=2\{R_{1212}^{~~~~2}+R_{1313}^{~~~~2}+R_{1414}^{~~~~2}+R_{2323}^{~~~~2}+R_{2424}^{~~~~2}+R_{3434}^{~~~~2}+R_{1212}R_{1313}\\
    &\qquad\quad+R_{1313}R_{1414}+R_{1212}R_{1414}+R_{1212}R_{2323}+R_{2323}R_{2424}+R_{1212}R_{2424}\\
    &\qquad\quad+R_{1313}R_{2323}+R_{2323}R_{3434}+R_{1313}R_{3434}+R_{1414}R_{2424}+R_{2424}R_{3434}\\
    &\qquad\quad+R_{1414}R_{3434}+\rho_{12}^{~~2}+\rho_{13}^{~~2}+\rho_{14}^{~~2}+\rho_{23}^{~~2}+\rho_{24}^{~~2}+\rho_{34}^{~~2}\},\\
    &\tau^2=4\{R_{1212}^{~~~~2}+R_{1313}^{~~~~2}+R_{1414}^{~~~~2}+R_{2323}^{~~~~2}+R_{2424}^{~~~~2}+R_{3434}^{~~~~2}+2R_{1212}R_{1313}\\
    &\qquad+2R_{1212}R_{1414}+2R_{1212}R_{2323}+2R_{1212}R_{2424}+2R_{1212}R_{3434}\\
    &\qquad+2R_{1313}R_{1414}+2R_{1313}R_{2323}+2R_{1313}R_{2424}+2R_{1313}R_{3434}\\
    &\qquad+2R_{2424}R_{3434}+2R_{1414}R_{2323}+2R_{1414}R_{2424}+2R_{1414}R_{3434}\\
    &\qquad+2R_{2323}R_{2424}+2R_{2323}R_{3434}\}.
    \eal
    \ee
From \eqref{eq:curv_1}, \eqref{eq:norm_curv} and \eqref{eq:Ric}, we
obtain the following equality
    \be\label{eq:id_{11}}
    \sum_{a,b,c}R_{abc1}^{~~~~2}-2\sum_{a}\rho_{1a}^{~~2}-2\sum_{a,b}\rho_{ab}R_{1ab1}+\tau\rho_{11}-\frac{1}{4}|R|^2+|\rho|^2-\frac{\tau^2}{4}=0.
    \ee
Similarly, we obtain
    \be\label{eq:id_{dd}}
    \sum_{a,b,c}R_{abcd}^{~~~~2}-2\sum_{a}\rho_{da}^{~~2}-2\sum_{a,b}\rho_{ab}R_{dabd}+\tau\rho_{dd}-\frac{1}{4}|R|^2+|\rho|^2-\frac{\tau^2}{4}=0
    \ee
 for $d=2$, $3$, $4$. Further, we get the following
    \be\label{eq:curv_12}
    \bal
    &\sum_{a,b,c}R_{abc1}R_{abc2}=2\{R_{1414}R_{1424}+R_{1423}R_{2324}+R_{1324}R_{2324}+R_{1424}R_{2424}\\
    &\qquad\qquad\qquad\qquad+R_{1334}R_{2334}+R_{1434}R_{2434}\}\\
    &\qquad\qquad\qquad\quad=2\{-\rho_{12}R_{1414}-\rho_{34}R_{1423}-\rho_{34}R_{1324}-\rho_{12}R_{2424}\\
    &\qquad\qquad\qquad\qquad+\rho_{14}\rho_{24}+\rho_{13}\rho_{23}\},\\
    &\sum_{a}\rho_{1a}\rho_{2a}=\rho_{11}\rho_{21}+\rho_{12}\rho_{22}+\rho_{13}\rho_{23}+\rho_{14}\rho_{24}\\
    &\qquad\qquad~=-\rho_{12}(2R_{1212}+R_{1313}+R_{1414}+R_{2323}+R_{2424})\\
    &\qquad\qquad\qquad+\rho_{13}\rho_{23}+\rho_{14}\rho_{24},\\
    &\sum_{a,b}\rho_{ab}R_{1ab2}=\rho_{12}R_{1212}-\rho_{34}(R_{1324}+R_{1423})+\rho_{44}\rho_{12}\\
    &\qquad\qquad\quad~=\rho_{12}(R_{1212}-R_{1414}-R_{2424}-R_{3434})-\rho_{34}(R_{1324}+R_{1423}),\\
    &\tau\rho_{12}=-2\rho_{12}(R_{1212}+R_{1313}+R_{1414}+R_{2323}+R_{2424}+R_{3434}).
    \eal
    \ee
Thus, from \eqref{eq:curv_12}, we have
    \be\label{eq:id_{12}}
    \sum_{a,b,c}R_{abc1}R_{abc2}-2\sum_{a}\rho_{1a}\rho_{2a}-2\sum_{a,b}\rho_{ab}R_{1ab2}+\tau\rho_{12}=0.
    \ee
Similarly, we see that
    \be\label{eq:id_{ij}}
    \sum_{a,b,c}R_{abci}R_{abcj}-2\sum_{a}\rho_{ia}\rho_{ja}-2\sum_{a,b}\rho_{ab}R_{iabj}+\tau\rho_{ij}=0
    \ee
holds for $(i,j)=(1,3)$, $(1,4)$, $(2,3)$, $(2,4)$, $(3,4)$. From
\eqref{eq:id_{11}}, \eqref{eq:id_{dd}}, \eqref{eq:id_{12}} and
\eqref{eq:id_{ij}}, we have \eqref{eq:id} with respect to a Chern
basis $\{e_i\}$. Since \eqref{eq:id} is a tensor equation,
\eqref{eq:id} is valid for any orthonormal basis. This completes the
proof of Main Theorem. \hfill$\qed$\medskip

As an application of our Main Theorem, we shall provide a new proof
of the following classical Theorem \cite{we}.
\begin{cor} Let
$M'=(M',g')$ be any 3-dimensional Riemannian manifold, and denote by
$R'$, $\rho'$ and $\tau'$ the curvature tensor, the Ricci tensor and
the scalar curvature of $M'$, respectively. Then from \eqref{eq:id},
we can see that the following equality
    \be\label{eq:3-dim-id}
    \bal
    R'_{abcd}=&\rho'_{ad}\delta_{bc}-\rho'_{ac}\delta_{bd}+\delta_{ad}\rho'_{bc}-\delta_{ac}\rho'_{bd}-\frac{\tau}{2}(\delta_{ad}\delta_{bc}-\delta_{ac}\delta_{bd})
    \eal
    \ee
$(1\leqq a,b,c,d\leqq 3)$ holds with respect to any orthonormal
basis $\{e'_a\}$ of $T_{p'}M'$ at any point $p'\in M'$.
\end{cor}
\begin{pf} Let $M$ be the Riemannian product of $M'$ and a real line
$\mathbb{R}$ and $\{e_i\}=\Big\{e_a=e'_a,e_4=\frac{d}{d t}\Big\}$ be
any orthonormal basis of $p=(p',t)\in M'\times \mathbb{R}$, where
$\{e'_a\}=\{e'_1,e'_2,e'_3\}$ is an orthonormal basis of $T_{p'}M'$.
Now, setting $i=j=4$ in \eqref{eq:id}, we can easily get the
following equality
    \be\label{eq:id_3_dim}
    \frac{1}{4}|R'|^2-|\rho'|^2+\frac{\tau'^2}{4}=0.
    \ee
Then, we can see that the equality \eqref{eq:id_3_dim} is equivalent
to the following
    \be\label{eq:3-dim-squre}
    \bal
    \sum_{a,b,c,d}\Big\{R'_{abcd}-\rho'_{ad}\delta_{bc}+\rho'_{ac}\delta_{bd}-\delta_{ad}\rho'_{bc}+\delta_{ac}\rho'_{bd}
    +\frac{\tau}{2}(\delta_{ad}\delta_{bc}-\delta_{ac}\delta_{bd})\Big\}^2=0.
    \eal
    \ee
Thus, from \eqref{eq:3-dim-squre}, we have the equality
\eqref{eq:3-dim-id}. \hfill$\qed$\medskip
\end{pf}

An $n$-dimensional Einstein manifold $M=(M,g)$ is called
super-Einstein if $M$ satisfies
    \be\label{eq:super-Einstein}
    {R^{abc}}_{i}R_{abcj}=\frac{1}{n}|R|^2g_{ij}
    \ee
\cite{GW}. We here remark that the constancy of $|R|^2$ follows from
the condition for an $n(\ne4)$-dimensional super-Einstein manifold
(\cite{B-V}, Lemma 3.3). For a 4-dimensional super-Einstein
manifold, the constancy of $|R|^2$ is usually required \cite{GW}.
Then from Main Theorem, we have immediately the following.
\begin{cor}\label{cor:4-dim einstein}
{{A 4-dimensional Einstein manifold satisfies the condition
\eqref{eq:super-Einstein}.}}
\end{cor}
We note that Corollary \ref{cor:4-dim einstein} can be also proved
by making use of a Singer-Thorpe basis in the 4-dimensional Einstein
manifold.  We shall call a Riemannian manifold $M=(M,g)$ satisfying
the condition \eqref{eq:super-Einstein}  a {\it weakly
super-Einstein manifold}. Especially, from the Corollary
\ref{cor:4-dim einstein}, we shall call a 4-dimensional Riemannian
manifold $M=(M,g)$ satisfying the condition
\eqref{eq:super-Einstein} (with $|R|^2$ not necessarily constant) a
{\it weakly Einstein manifold} in short. The following example shows
that a weakly Einstein manifold is not necessarily Einstein.
\begin{ex}
Let $M$ be a Riemannian product manifold of 2-dimensional Riemannian
manifolds $M_1^2(c)$ and $M_2^2(-c)$ of constant Gaussian curvatures
$c$ and $-c$ $(c\ne0)$, respectively. Then we can easily check that
$M$ is not Einstein.
 We also can easily check that $M$ satisfies \eqref{eq:super-Einstein}, thus
 $M$ is weakly Einstein.
 \end{ex}

The detailed study of the weakly Einstein spaces is in procedure and
it will be published elsewhere.

\section*{Acknowledgements}

 Research of Yunhee Euh was supported by the National
Research Foundation of Korea Grant funded by the Korean Government
[NRF-2009-352-C00007]. Research of JeongHyeong Park was supported by
Basic Science Research Program through the National Research
Foundation of Korea(NRF) funded by the Ministry of Education,
Science and Technology (2009-0087201).


\end{document}